\documentclass[11pt,twoside]{article}

\usepackage{a4wide}
\usepackage{amssymb}
\usepackage[matrix,arrow,curve]{xy}
\usepackage{epsfig}

\usepackage{amsmath}
\usepackage{color}
\usepackage{hyperref}
 \definecolor{darkblue}{rgb}{0,0,.7}
 \hypersetup{pdftex=true, colorlinks=true, breaklinks=true, linkcolor=darkblue, citecolor=darkblue, menucolor=darkblue, pagecolor=blue, urlcolor=darkblue}

\numberwithin{equation}{section}
\numberwithin{figure}{section}

\usepackage{theorem}
\newtheorem{thm}{Theorem}[section]
\newtheorem{lem}[thm]{Lemma}
\newtheorem{pro}[thm]{Proposition}
\newtheorem{cor}[thm]{Corollary}

\newtheorem{remark}[thm]{Remark}

\newtheorem{defi}[thm]{Definition}

\newcommand{\comment}[1]{}

\newcommand{\Z}{\mathbb Z}

\newcommand{\Pro}{\mathrm{P}^{\infty}}
\newcommand{\Po}{{\mathrm{P}}}
\newcommand{\HB}{\mathbb H}

\newcommand{\R}{\mathbb R}

\newcommand{\C}{\mathbb C}
\newcommand{\qed}{\quad \textbf{q.e.d.}}

\DeclareMathOperator{\Tor}{Tor}

\xyoption{frame}
\xyoption{tips}
\xyoption{color}
\entrymodifiers={+!!<0pt,\fontdimen22\textfont2>}

\title{Involutions on $S^6$ with 3-dimensional fixed point set}
\author{Martin Olbermann}
\date{\today}
\begin{document}
\maketitle

\begin{abstract}
\noindent In this article, we classify all involutions on $S^6$ with 3-dimensional fixed point set.
In particular, we discuss the relation between the classification of involutions with fixed point set a knotted 3-sphere
and the classification of free involutions on homotopy $\C \Po^3$'s.
\end{abstract}


\section{Introduction}
As a general assumption, we are interested in smooth involutions on connected, closed, smooth manifolds.

The study of group actions on very simple manifolds such as disks, spheres or Euclidean spaces
has been a very active subject. In his MathSciNet review of \cite{Paw}, Masuda notes: 
``Representations of groups provide examples of group actions on Euclidean spaces, disks or spheres, 
and an important natural question in transformation groups is to what extent arbitrary actions on those 
spaces resemble actions provided by representations."
The first highlight theorems are due to P.A. Smith. For an involution on a sphere, 
Smith proved that the fixed point set is a $\Z_2$-homology sphere. 
There are various converses to this theorem which can be found in the literature, 
e.g. \cite{DW, Paw} and references therein. However, the following theorem seems to be new:
\begin{thm} \cite{confree}
Every $\Z_2$-homology 3-sphere is the fixed point set of an involution on $S^6$.
\end{thm}
Moreover, the method used in \cite{confree} generalizes to a classification of these involutions, 
which is the subject of the present paper. (A shorter proof of the above existence theorem can be
given using Dovermann's equivariant surgery approach along the lines of \cite{Sch} and \cite{DMS}.)

\bigskip

Recently, a class of involutions called {\em{conjugations}} was defined in \cite{HHP} and
various aspects of conjugations were studied in \cite{FP, OlbThesis, HHo, HH}. 
Conjugations $\tau$ on topological spaces $X$ have the property that the
fixed point set has $\Z_2$-cohomology ring isomorphic to the $\Z_2$-cohomology ring of $X$, 
with the slight difference that all degrees are divided by two. 
In the case of smooth involutions on (positive-dimensional) spheres, the conjugations are exactly the involutions on even-dimensional
spheres with half-dimensional fixed point set. 
In dimension 2, the Sch\"onflies theorem gives a classification: every conjugation is conjugate to the reflection of $S^2$ at the equator.
In dimension 4, work of Gordon and Sumners shows that there are infinitely many non-equivalent 
conjugations on $S^4$. Hambleton and Hausmann recently reduced the study of such involutions to a non-equivariant
four-dimensional knot theory question \cite{HH}. 
Knot theory of $k$-spheres in $S^{2k}$ is easier for $k>2$, so that the study of
conjugations gets a different flavor for $k>2$.
 
\bigskip

The study of free involutions on simply-connected spin manifolds with the same homology groups as $\C \Po^3$ 
was motivated by the question whether it is possible to define an ``equivariant Montgomery-Yang correspondence". 
After Haefliger \cite{Hae} proved that the group $C^3_3$ of knotted $S^3$'s in $S^6$ 
(isotopy classes of smooth embeddings, or equivalently
diffeomorphism classes relative to $S^3$) is isomorphic to $\Z$,
Montgomery and Yang showed \cite{MY} that there is a natural bijection between $C^3_3$ and the set of 
diffeomorphism classes of homotopy $\C \Po^3$'s.  
Wall's classification of all simply-connected spin manifolds with the same homology groups as $\C \Po^3$ \cite{Wal} also uses the bijection between 
diffeomorphism classes of such manifolds, together with a basis of $H^2$, and isotopy classes of framed $S^3$-knots in $S^6$.
In the equivariant setting, Li and L\"u \cite{LL} show that the existence of a free involution on a homotopy 
$\C \Po^3$ implies the existence of an involution on $S^6$ with fixed point set the corresponding knotted $S^3$.
Similarly, in our approach, the same surgery arguments apply in both cases.
However, it is not possible to produce a nice bijection on the set of equivariant
diffeomorphism classes of these, as claimed in \cite{Su}. Our classification results are:

\begin{thm}\label{Mthm}
Let $M$ be a smooth closed simply-connected spin manifold
with $H_2(M)=\Z$, and $H_3(M)=0$. Let $x\in H^2(M;\Z)$ be a generator. 
We assume that $M\not \cong S^2\times S^4$.
\begin{itemize}
\item
If $\left\langle \frac{p_1(M)x-4x^3}{24},[M]\right\rangle$
is odd, there exists no free involution on $M$.
\item
If $\langle x^3,[M]\rangle$ is odd, and $\left\langle \frac{p_1(M)x-4x^3}{24},[M]\right\rangle$
is even, there exist up to diffeomorphism exactly two free involutions on $M$.
\item
If $\langle x^3,[M]\rangle$ is even, and $\left\langle \frac{p_1(M)x-4x^3}{24},[M]\right\rangle$
is even, there exist up to diffeomorphism exactly five free involutions on $M$.
\end{itemize}
If $M\cong S^2\times S^4$, then the same classification holds for orientation-reversing involutions
which act by $-1$ on $H^2(M)$.
\end{thm}
In the case of homotopy $\C P^3$'s the classification of free involutions 
was given by Petrie \cite{Pet} (whose result contains a mistake, corrected by Dovermann, Masuda and Schultz), and Su \cite{Su}.
Li and Su (unpublished) give the answer to the existence question in the case of odd $\langle x^3,[M]\rangle$. 
Our method reproves all these results in a different way,
extends to a larger class of manifolds, and gives classification results in all cases.

\begin{thm}\label{thm13}
For every even element of $C^3_3$ there are (up to equivariant diffeomorphism relative $S^3$) exactly four involutions (conjugations)
on $S^6$ with the knot as fixed point set. For every odd element of $C^3_3$, there is no involution on $S^6$ with the knot 
as fixed point set.
\end{thm}
The new part of the theorem is the classification of these involutions. Li and L\"u proved that the existence of an involution
with fixed point set a given knot is equivalent to the existence of a free involution on the corresponding homotopy $\C \Po^3$
under the Montgomery-Yang correspondence, so that together with Su's work mentioned above, the existence question was answered.

\bigskip

The case of involutions on $S^6$ with fixed point set $S^3$ is especially interesting since, given another 6-manifold with an involution that has a
3-dimensional fixed point set, (equivariant) connected sum gives a possibly different involution on the 
same 6-manifold with same fixed point set.
(This is in analogy with the fact that connected sum with a homotopy sphere gives a possibly
new smooth structure on the same underlying topological manifold.)
Similarly, connected sum with a conjugation on $S^6$ with fixed point set different from $S^3$
gives a new conjugation on the same 6-manifold, with different fixed point space.

\bigskip

Our main result is a classification of smooth involutions on $S^6$ with arbitrary three-dimensional fixed point set, 
using a recent classification of embeddings of closed oriented connected 3-manifolds into $S^6$ by A.~Skopenkov \cite{Sko}.
Skopenkov proves that for a 3-dimensional $\Z_2$-homology sphere $M$ the set of isotopy classes of embeddings $i:M\to S^6$ 
has a free action by $C_3^3$ and the orbits are in canonical bijection with $H_1(M)$.

\begin{thm}\label{embcl} Let $M$ be a $\Z_2$-homology sphere of dimension 3.
The set of isotopy classes of embeddings $i:M\to S^6$ which are the fixed point sets of involutions (conjugations) is contained in the orbit
corresponding to $0\in H_1(M)$. Moreover, it is acted upon freely and transitively by 
$2\Z\subseteq \Z\cong C_3^3$. There are up to equivariant diffeomorphism relative to $i$ exactly four such conjugations for every such $i$. 
\end{thm}

In our case, we can also classify involutions without additional identification of the fixed point set with a given 3-manifold.
However, we consider involutions together with an orientation of their fixed point sets,
and the equivalence relation is equivariant diffeomorphism which respects the orientations of both $S^6$ and the fixed point set. 
Equipping the involution with an orientation of the fixed point set is necessary in order to perform 
connected sums. Thus it is more natural to determine this more structured set of equivariant diffeomorphism classes
$Inv_{M}(S^6)$. 

Since the action of the mapping class group of $M$ on the above set of equivariant diffeomorphism classes of involutions relative to the fixed
point set $M$ is trivial, we get the same classification as in theorem \ref{embcl}: 
\begin{thm}\label{invcl}
Conjugations up to conjugation (preserving orientations) with fixed point set of
a fixed oriented diffeomorphism type of $\Z_2$-homology 3-spheres are in bijection with $\Z \oplus \Z_4$.
Under connected sum, $Inv_{S^3}(S^6)\cong \Z \oplus \Z_4$ is an isomorphism of groups, and
$Inv_{S^3}(S^6)$ acts freely and transitively on $Inv_{M}(S^6)$ for every $\Z_2$-homology 3-sphere $M$.
\end{thm}

\begin{remark}
In principle, using the machinery described in \cite{KrM} it is also possible to prove 
classification results for conjugations on $X^6$ with fixed point set $M^3$
in other cases as $X=S^6$. 
(However, the argument we use to show the surgery obstruction is trivial 
does not extend to the case of free involutions on other manifolds.)

One would compute the set of equivariant diffeomorphism classes of the involution together with an identification of the fixed point set with a 
prescribed 3-manifold $M$, i.e. the set
$$ Emb_{\Z_2}(M\to X)=\{ f:(M,id)\to (X,\tau) \} / \sim $$
where $f$ is an inclusion onto the fixed point set of $\tau$, and  $f:(M,id)\to (X,\tau)$
is identified with $f':(M,id)\to (X',\tau')$ if there is an equivariant diffeomorphism $\phi:(X,\tau)\to
(X',\tau')$ making the obvious triangle commute.
One of the difficulties to overcome is ``due" to Wall's classification: in order to determine which of the resulting 6-manifolds
are diffeomorphic to $X$, one would need a good understanding of the isomorphism classes of the algebraic invariants (trilinear forms), 
and this problem seems to be very hard in general.
\end{remark}

\noindent {\bf{Acknowledgements.}} I would like to thank Diarmuid Crowley, Jean-Claude Hausmann, Matthias Kreck and Arturo Prat-Waldron  
for many helpful discussions and remarks. Special thanks to Yang Su, whose article \cite{Su} was the origin of this paper, and who shared
with me the so far unpublished notes on a generalization of \cite{Su} due to himself and Banghi Li.

\section{Preliminaries}
\subsection{Modified surgery}

We will use Kreck's modified surgery theory \cite{Kre, KrM} which is also suited to give classification results.
By the equivariant tubular neighbourhood theorem, we can write $S^6=M\times D^3 \cup_\partial V$, where 
$V$ is a manifold with boundary $M\times S^2$ equipped with a free involution, which restricts to $(id,-id)$ on the boundary.
(It is not hard to see that the normal bundle of $M$ in $S^6$ is trivial. See \cite{OlbThesis} for a proof.)

A classification of manifolds $V$ with free involutions $\tau$ up to equivariant diffeomorphism is the same as a classification of the quotient manifolds 
$W=V/\tau$ up to diffeomorphism. This is what modified surgery theory will give us.

We first determine the normal 2-type of the manifolds $W$ under consideration.
The normal 2-type $B$ of a 6-manifold (see the precise definition \ref{normal2type})
is a fibration $B\to BO$. It
carries roughly the information of a 3-skeleton of the manifold together with the restriction
of the normal bundle to this 3-skeleton. 
After computing the bordism group of manifolds with normal $B$-structure, we show that in every bordism class there exists
a manifold (together with a map to $B$) which qualifies as the $W$ above.
Moreover, we show that given two normally $B$-bordant manifolds $W$ as above,  
the obstruction, which a priori lies in the complicated monoid $l_7(\Z_2,-1)$,
for the existence of an s-cobordism (i.e. a diffeomorphism) is zero. 
The diffeomorphism classification of the manifolds $W$ under consideration is given by the set of orbits of the action of 
the group of fiber homotopy self-equivalences $B\to BO$.

\subsection{Conjugations on manifolds}

This section explains what conjugation spaces are and shows that the smooth involutions on $S^{2n}$ with $n$-dimensional fixed point set 
are exactly the smooth conjugations. The rest of the paper does not depend on the material in this section.

\bigskip

A conjugation on a topological space $X$ is an involution $\tau:X\to X$, which we consider as an action of the group $\Z_2\cong C= \{id, \tau\}$
on $X$, and which satisfies the following cohomological pattern:
We denote the Borel equivariant cohomology of $X$ by $H^*_C(X;\Z_2)$. It is a module over $H^*_C(pt;\Z_2)=\Z_2[u]$. 
The restriction maps in equivariant cohomology are denoted by $\rho: H^*_C(X;\Z_2)\to H^*(X;\Z_2)$
and $r: H^*_C(X;\Z_2)\to H^*_C(X^\tau;\Z_2)\cong H^*(X^\tau;\Z_2)[u]$.

\begin{defi}\cite{HHP}\label{conjsp}
$X$ is a conjugation space if 
\begin{itemize}
\item
$H^{odd}(X;\Z_2)=0$, 
\item there exists a (ring) isomorphism $\kappa: H^{2*}(X;\Z_2)\to H^*(X^\tau;\Z_2)$
\item and a (multiplicative) section $\sigma: H^*(X;\Z_2)\to H^*_C(X;\Z_2)$ of $\rho$ 
\item such that the so-called conjugation equation holds:
$$r\sigma(x)=\kappa(x)u^k + \text{ terms of lower degree in }u.$$
\end{itemize}
\end{defi}
One does not need to require that $\kappa$ and $\sigma$ be ring homomorphisms, it is a consequence of the definition. Moreover,
the ``structure maps" $\kappa$ and $\sigma$ are unique, and natural with respect to equivariant maps between conjugation spaces.

There are many examples of such conjugations: complex conjugation on the projective space $\C \Po^n$ and on
complex Grassmannians, natural involutions on smooth toric manifolds \cite{DJ} and on 
polygon spaces \cite{HK}. Every cell complex with the property that each cell is a unit disk in $\C^n$ with complex conjugation, 
and with equivariant attaching maps, is a conjugation space. 
Coadjoint orbits of semi-simple Lie groups with the Chevalley involution are conjugation spaces.
Moreover, there are various constructions of new conjugation spaces out of other conjugation spaces.

A conjugation manifold is a conjugation space consisting of a smooth manifold $X$ with a smooth involution $\tau$.
As a consequence, a closed conjugation manifold $X$ must be even-dimensional, say of dimension $2n$, and $M$ is of dimension $n$. 

In \cite{OlbThesis} we proved that it is possible to give a definition of conjugation spaces without
the non-geometric maps $\kappa$ and $\sigma$, which is moreover well-adapted to the case of conjugation manifolds,
where the fixed point set has an equivariant tubular neighbourhood. 

\begin{pro}
Every smooth involution on $S^{2n}$ with (non-empty) n-dimensional fixed point set is a conjugation.
\end{pro}
Proof:  Let $pt\in S^{2n}$ be a fixed point of the involution.
Then $pt$ is a conjugation space and $(S^{2n},pt)$ is a conjugation pair, by 
the same proof as in Example 3.5 of \cite{HHP}. Then the extension property for triples, 
Prop. 4.1 in \cite{HHP}, shows that $S^{2n}$ is a conjugation space. 
(A slightly different proof is given in \cite{HH}.)
\qed

\section{Free involutions on certain 6-manifolds}

We are considering smooth involutions on smooth closed simply-connected spin manifolds $M$
with $H_2(M)=\Z$, and $H_3(M)=0$. The classification by Wall in \cite{Wal}
also uses a generator $x\in H^2(M)$  and an orientation of $M$.
Then pairs $(M,x)$ up to diffeomorphism are classified by the bordism class
of $M\stackrel x\to \C \Pro \in \Omega_6^{Spin}(\C\Pro)\cong \Z^2$,
and $(M,x)$ is mapped under the isomorphism to 
$$\left(  \left\langle \frac{p_1(M)x-4x^3}{24},[M]\right\rangle , \langle x^3,[M]\rangle\right).$$
Both switching the sign of the generator of $H^2(M)$ and the orientation of $M$ 
induce multiplication with $-1$, so that diffeomorphism classes are in bijection with 
$\Z^2/-1$.

As observed in \cite{LL}, the Lefschetz fixed point theorem implies:
\begin{lem}
If $\langle x^3,[M]\rangle$ is non-zero, then a free involution on $M$ must be 
orientation reversing, and act by $-1$ on $H^2(M)$.
\end{lem}

If $\langle x^3,[M]\rangle=0$, we consider only orientation reversing free involutions
which are $-1$ on $H^2(M)$. From Li and Su we learned that except for the case 
$M=S^2\times S^4$, these are all free involutions: If the involution acts by $-1$ on $H^4(M)$,
the first Pontryagin class must be 0, and we use the classification. 
The remaining case is handled as above by the Lefschetz fixed point theorem.

Obviously, for $M=S^2\times S^4$, our classification of orientation reversing free involutions
which are $-1$ on $H^2(M)$ does not include all free involutions.

\subsection{The normal 2-type}
\begin{defi}\label{normal2type}
The normal 2-type of a compact manifold $N$ is a fibration $B_2(N)\to BO$ which is obtained
as a Postnikov factorization of the stable normal bundle map $N\to BO$:
There is a 3-connected map $N\to B_2(N)$, the fibration 
$B_2(N)\to BO$ is 3-coconnected (i.e. $\pi_i(BO,B_2(N))=0$ for $i>3$), and
the composition is the stable normal bundle map. This determines $B_2(N)\to BO$ up to fiber homotopy equivalence.
\end{defi}

\begin{lem}
Let $\tau$ be an involution on $M$ as above, 
and let $N=M/\tau$ be the the quotient space of the involution.
The second space in a Postnikov tower for $N$ is 
either $P=(\C\Pro \times S^\infty)/(c,-1)$, where
$c$ is complex conjugation, or $Q=(\C\Pro \times S^\infty)/(\tau,-1)$, where
$\tau$ is fiberwise the antipodal involution on $S^2\to \C \Pro \to \HB \Pro$.
\end{lem}
Proof: The first space in the Postnikov tower is a $K(\Z_2,1)$, and the second space is a $K(\Z,2)$ fibration 
over it, with $\pi_1$ acting nontrivially on $\pi_2$. Such fibrations are classified by their $k$-invariant in 
$H^3(K(\Z_2,1);\Z_-)\cong \Z_2$. The spaces $P$ and $Q$ have the required properties, and they 
are not homotopy equivalent, as e.g. $H^2(P;\Z_2)\cong \Z_2^2$ and $H^2(Q;\Z_2)\cong \Z_2$. 
So they represent all isomorphism classes of fibrations. ($P$ has $k$-invariant 0, and $Q$ has nonzero $k$-invariant.)
\qed

\begin{lem}\label{cohQ}
The $\Z_2$ cohomology ring of $N$ is $\Z_2[q,t]/\langle t^3,q^2\rangle$, where $deg(q)=4, deg(t)=1$. 
\end{lem}
Proof: 
(This was proved in \cite{Su} in the case of homotopy $\C \Po^3$'s, and we generalize his proof.)
We consider the Serre spectral sequence of the fibration $M\to N \to \R\Pro$, with $\Z_2$ (and also sometimes
with integral) coefficients. The first case is that $d_3:E_3^{0,2}\to E_3^{3,0}$ is non-trivial.
Then by multiplicativity the $E_4$-term has exactly one $\Z_2$ in $E_4^{p,0}$ and $E_4^{p,4}$ for each $p=0,1,2$.
There are no further differentials, and we get the above cohomology ring.
The second case is that $d_3:E_3^{0,2}\to E_3^{3,0}$ is trivial. We will show that this leads to a contradiction.
By multiplicativity, also $d_3:E_3^{0,6}\to E_3^{3,4}$ is trivial. If we remember that we need the limit to have no
cohomology in degrees $>6$, then we see in the sequence with integral coefficients that there must be a nontrivial
$d_3$-differential between the fourth and second line. Then the same must hold for the $\Z_2$ coefficient
sequence. And we get a $d_7$-differential from the sixth to the zeroth line.
The $E_\infty$-term has exactly one $\Z_2$ in $E_\infty^{p,0}$ for $p=0, \dots , 6$ and $E_\infty^{p,2}$ for $p=0,1,2$.

This gives a cohomology ring with a generator $t\in H^1(N;\Z_2)$,
and another generator $x\in H^2(N;\Z_2)$. Since in this case $H^2(N;\Z_2)\cong \Z_2^2$, the second
Postnikov space must be $P$, and we can choose $x\in H^2(N;\Z_2)$ coming from $P$; 
we choose the class in $H^2(P;\Z_2)$ which maps to 0 under a section of $P\to \R \Pro$. 
Note that $x$ maps nontrivially to $H^2(M;\Z_2)$.

We have the relations $t^7=t^3x=x^2+at^2x+t^4=0$, where $a\in\Z_2$. 
(We have $Sq^1x=tx$ as this is true in $P$, 
thus $t^3x=Sq^1(t^2x)=0$ since $H^5(N;\Z)=0$.
By Poincar\'e duality $t^2x^2$ can't be zero. This implies that in $x^2+at^2x+bt^4=0$ we have $b=1$.)

We see that $Sq^1(t^5)=t^6, Sq^2(t^4)=0, Sq^2(t^2x)=t^2x^2\ne 0$. It follows that the first Wu class is $t$ 
and the second Wu class is $x$. But that implies that the second Stiefel-Whitney class of $N$ maps
to a non-trivial class in $H^2(M;\Z_2)$. But since $M$ is spin, this image must be zero, as it is the second
Stiefel-Whitney class of $M$. Contradiction. \qed

\begin{cor}
The second space in a Postnikov tower for $N$ is $Q$.
\end{cor}
Proof: This follows from the fact that $H^2(N;\Z_2)\cong \Z_2$. \qed

\begin{pro}
If such a manifold $M$ with $\langle x^3,[M]\rangle$ odd has
a free involution $\tau$, then the quotient space $N$ has normal 2-type 
$B =  BSpin \times Q \to BO \times BO(1)\stackrel{\oplus}\to BO$. 
\end{pro}
Proof: If $\langle x^3,[M]\rangle$ is odd, then the map
$H^4(Q;\Z_2)\to H^4(N;\Z_2)$ is a bijection, since then we have isomorphisms
$H^4(Q;\Z_2)\to H^4(\C \Pro;\Z_2)\to H^4(M;\Z_2) \leftarrow H^4(N;\Z_2)$.
Then $Sq^2:H^4(N;\Z_2)\to H^6(N;\Z_2)$ is zero, so is the second Wu class of
$M/\tau$, and as a consequence the quotient space has a spin structure twisted by $L\to Q$.
\qed

\begin{pro}\label{normalbpro}
If such a manifold $M$ with trivial square $H^2(M;\Z_2)\to H^4(M;\Z_2)$ has
a free involution $\tau$, then the quotient space $N$ has one of the following normal 2-types:
\begin{eqnarray*}
A & = & Q\times BSpin \to BO(1)\times BO(1)\times BO(1)\times BO\stackrel{\oplus}\to BO,\\
B & = & Q\times BSpin \to BO(1)\times BO\stackrel{\oplus}\to BO,\\
\end{eqnarray*}
where the maps to all $BO(1)$'s are the projections $p:Q\to \R\Pro=BO(1)$.
\end{pro}
Proof: If $\langle x^3,[M]\rangle$ is even, then there is a second possibility for the normal 2-type.
If we fix the second space in a Postnikov tower to be $Q$, 
the second Stiefel-Whitney class of $N$ can be $t^2$ or 0. Thus $N$ either
admits spin structures twisted by $L$ or spin structures twisted by $L\oplus L \oplus L = 3L$.
\qed

\bigskip 

\noindent
\label{normalBstr}
A normal $B$-structure on a manifold $N$ can be defined in three equivalent ways. 
\begin{itemize}
\item
It is a vertical homotopy class of lifts of the normal bundle map $N\to BO$ to $B$
(this is independent of the particular map $N\to BO$ coming from an embedding of $N$ into some Euclidean space).
\item
It is a map $f:N\to Q$ (up to homotopy) together with a spin structure on the bundle $\nu_N - f^*(L)$, where $L$ is the nontrivial 
real line bundle on $Q$.
\item
It is a map $f:N\to Q$ together with a homotopy $\eta$ (and this up to homotopy) in the following square:
$$ 
\xymatrix{
N
\ar[r]^\nu
\ar[d]^f
& BO
\ar[d]^{w_1\times w_2}
\\
Q
\ar@{=>}[ur]^{\eta}
\ar[r]_/-2em/{t \times 0}
& K(\Z_2,1)\times K(\Z_2,2)
}
$$
Here we fix maps corresponding to the classes 
$w_1\in H^1(BO;\Z_2), w_2\in H^2(BO;\Z_2)$, the generator $t\in H^1(Q;\Z_2), 0\in H^2(Q;\Z_2)$.

(For fixed $f$, the homotopy classes of homotopies $\eta$ have a free and transitive action by 
$\pi_1((K(\Z_2,1)\times K(\Z_2,2))^N) \cong H^1(N;\Z_2)\times H^0(N;\Z_2)$.) 
\end{itemize} 
Similarly, normal $A$-structures are defined. 
(In the first definition, replace $B$ by $A$. In the second definition, replace $L$ by $3L$. In the third
definition, replace $0$ by $t^2$.)

\begin{remark}
As a converse to proposition \ref{normalbpro}, a manifold with normal $B$-structure $N^6\to B$ which is a 3-connected map
is the quotient of an involution on a closed simply-connected spin manifold $M$
with $H_2(M)=\Z$, and $H_3(M)=0$ if and only if $H_3(N;\Lambda)=0$. Here $\Lambda=\Z[\Z_2]$ is the group ring of the fundamental group.
The same holds for normal $A$-structures.
\end{remark}

\subsection{Computation of the bordism groups}

We compute bordism groups of manifolds with normal $A$-structures (resp. $B$-structures).
The (co)homology of $Q$ is described in \cite{OlbThesis}. 
To compute the bordism groups $\Omega_6^A$ and $\Omega_6^B$ we use the 
Atiyah-Hirzebruch spectral sequence (which computes the group up to an extension problem) and
an Adams spectral sequence (which can help solve the extension problem).
We get:
\begin{thm}
We have isomorphisms $\Omega_6^A \cong \Z^2 \oplus \Z_2$ and 
$\Omega_6^B \cong \Z^2 \oplus \Z_4.$
\end{thm}
Proof: The case of $\Omega_6^B$ was proven in \cite{OlbThesis}.
For $\Omega_6^A$, the Atiyah-Hirzebruch spectral sequence is 
$$H_p(Q;\underline{\Omega_q^{Spin}})\Rightarrow \Omega_{p+q}^A,$$ 
the $d^2$-differential 
$$E^2_{p,1}\cong H_p(Q;\Z_2)\to H_{p-2}(Q;\Z_2)\cong E^2_{p-2,2}$$ is the dual of $Sq^2+tSq^1+t^2Sq^0$, 
and the $d^2$-differential 
$$E^2_{p,0}\cong H_p(Q;\Z_-)\to H_{p-2}(Q;\Z_2)\cong E^2_{p-2,1}$$ is reduction modulo 2 composed with 
the dual of $Sq^2+tSq^1+t^2Sq^0$.
(We obtain this using the Thom isomorphism $$\Omega_6^{Spin}(Q;3L)\cong \Omega_9^{Spin}(D(3L),S(3L))$$
and the Atiyah-Hirzebruch spectral sequence for the latter. See also \cite{Tei}.)
From the calculations in \cite{OlbThesis} it immediately follows that the differentials
$d^2:E^2_{6,1}\to E^2_{4,2}$ and $d^2:E^2_{6,0}\to E^2_{4,1}$ are non-trivial.
This implies that the nonzero terms on the sixth diagonal in the $E^\infty$-term are:
$$E^{\infty}_{2,4}\cong \Z,\quad E^{\infty}_{5,1}\cong \Z_2,\quad E^{\infty}_{6,0}\cong 2\Z.$$
Thus $\Omega_6^A  \cong  \Z^2 \oplus \Z_2$ or $\Omega_6^A  \cong  \Z^2$.

Now we consider the Adams spectral sequence
$$Ext_\mathcal{A}^{s,t}(H^*(MSpin\wedge T(3L);\Z_2),\Z_2)\Rightarrow \Omega_{t-s-3}^A,$$ 
where $\mathcal{A}$ is the mod 2 Steenrod algebra.
We compute the left hand side for $t-s-3\le 6$.
We have 
$$H^*(MSpin\wedge T(3L);\Z_2)\cong H^*(MSpin;\Z_2)\otimes \tilde{H}^*(T(3L);\Z_2),$$
and $\tilde{H}^*(T(3L);\Z_2)$ is a free $H^*(Q;\Z_2)$-module on one generator $u_3$ of degree 3 (the Thom class).
We have $$Sq(u_3)= w(3L)u_3 = u_3+tu_3+t^2u_3+t^3u_3.$$
All of this allows to write down the $\mathcal{A}$-module structure of $H^*(MSpin\wedge T(3L);\Z_2)$ in degrees $\le 10$.
Then we compute a free $\mathcal{A}$-resolution (in low degrees). From this we get the $E^2$-term of the spectral sequence, which is
displayed in the following diagram:
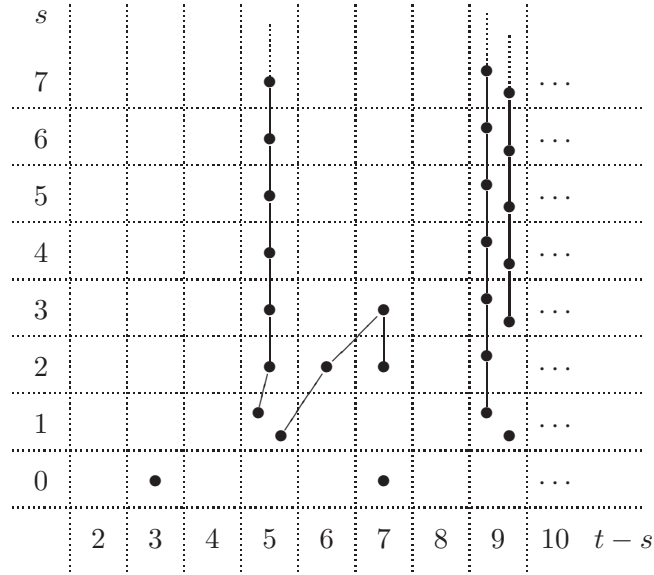
\[
 \xy <1.8pc,0pc>:<0pc,1.8pc>::
 (0.2,-0.05) ="A"; (-0.2,0.35)="B",(0.0,0.15) ="C",
 (-1,0) *!{0},
 (-1,1) *!{1},
 (-1,2) *!{2},
 (-1,3) *!{3},
 (-1,4) *!{4},
 (-1,5) *!{5},
 (-1,6) *!{6},
 (-1,7) *!{7},
 (-1,8.2) *!{s},
 (0,-1) *!{2},
 (1,-1) *!{3}, 
 (2,-1) *!{4},
 (3,-1) *!{5},
 (4,-1) *!{6},
 (5,-1) *!{7},
 (6,-1) *!{8},
 (7,-1) *!{9},
 (8,-1) *!{10},
 (9.2,-1) *!{t-s},
 (8,0) *!{\cdots},
 (8,1) *!{\cdots},
 (8,2) *!{\cdots},
 (8,3) *!{\cdots},
 (8,4) *!{\cdots},
 (8,5) *!{\cdots},
 (8,6) *!{\cdots},
 (8,7) *!{\cdots},
 (1,0)+"C" *{\bullet},
 (3,8)+"C"; (3,7)+"C" *{\bullet} **@{.}; (3,6)+"C" *{\bullet} **@{-};(3,5)+"C"="id42" *{\bullet} **@{-};
 (3,4)+"C"*{\bullet} **@{-}; (3,3)+"C" *{\bullet} **@{-}; (3,2)+"C" *{\bullet} **@{-}; (3,1)+"B" *{\bullet} **@{-},
 (3,1)+"A" *{\bullet}; (4,2)+"C"="id61" *{\bullet} **@{-}; (5,3)+"C" *{\bullet} **@{-};
 (5,2)+"C" *{\bullet} **@{-},
 (5,0)+"C" *{\bullet}, 
 (7,8)+"A"; (7,7)+"A" *{\bullet} **@{.}; (7,6)+"A" *{\bullet} **@{-}; (7,5)+"A" *{\bullet} **@{-};
 (7,4)+"A" *{\bullet} **@{-};(7,3)+"A" *{\bullet} **@{-},
 (7,8)+"B"; (7,7)+"B" *{\bullet} **@{.}; (7,6)+"B" *{\bullet} **@{-}; (7,5)+"B" *{\bullet} **@{-};
 (7,4)+"B" *{\bullet} **@{-};(7,3)+"B" *{\bullet} **@{-};(7,2)+"B" *{\bullet} **@{-};
 (7,1)+"B" *{\bullet} **@{-},
 (7,1)+"A" *{\bullet},
 (-1.5,-0.3) ; (9.5,-0.3) **@{.},
 (-1.5,0.7) ; (9.5,0.7) **@{.},
 (-1.5,1.7) ; (9.5,1.7) **@{.},
 (-1.5,2.7) ; (9.5,2.7) **@{.},
 (-1.5,3.7) ; (9.5,3.7) **@{.},
 (-1.5,4.7) ; (9.5,4.7) **@{.},
 (-1.5,5.7) ; (9.5,5.7) **@{.},
 (-1.5,6.7) ; (9.5,6.7) **@{.},
 (-0.5,-1.5) ; (-0.5,8.5) **@{.},
 (0.5,-1.5) ; (0.5,8.5) **@{.},
 (1.5,-1.5) ; (1.5,8.5) **@{.},
 (2.5,-1.5) ; (2.5,8.5) **@{.},
 (3.5,-1.5) ; (3.5,8.5) **@{.},
 (4.5,-1.5) ; (4.5,8.5) **@{.},
 (5.5,-1.5) ; (5.5,8.5) **@{.},
 (6.5,-1.5) ; (6.5,8.5) **@{.},
 (7.5,-1.5) ; (7.5,8.5) **@{.},
 \endxy
\]
Since there are no differentials starting or ending at $(t-s,s)=(9,1)$, we obtain
$\Omega_6^A  \cong  \Z^2 \oplus \Z_2$.
\qed
 
\subsection{Construction and classification up to normal $B$-bordism} \label{Con}

By Wall's classification, every bordism class in $\Omega_6^{Spin}(K(\Z,2))\cong \Z^2$ contains a unique 
normal 2-smoothing with $H_3=0$ (up to diffeomorphism relative $BSpin\times K(\Z,2)$).

\begin{thm}\label{freebord}
In every bordism class in $\Omega_6^A$ there is a unique manifold 
$W\to A$ (up to diffeomorphism relative to $A$) such that $W\to A$ is 3-connected
and $H_3(W;\Lambda)=0$.
The same is true if one replaces $A$ by $B$.
\end{thm}
Proof: For the construction and classification of free involutions on these manifolds, we use surgery theory.
The existence proof is a simplified version of the proof of the main theorem in \cite{confree}. 

We start with any 6-dimensional closed manifold with normal $A$-structure, 
and we do surgery to get manifolds $W$ such that the map $W\to A$ is 3-connected and $H_3(W;\Lambda)=0$.

Surgery below the middle dimension is always possible \cite{Kre}. It allows to modify any closed 6-manifold $W$ with
normal $A$-structure into one with normal 2-type $A$. Let us denote this new
manifold again by $W$. It remains to kill $H_3(W;\Lambda)$.

By the Hurewicz theorem (its extended version) we have a surjection
$\pi_3(W)\to H_3(W;\Lambda)$, where $\Lambda=\Z[\Z_2]$ is the group ring 
of the fundamental group, and $H_3(W;\Lambda)$ can be identified with the homology of the universal cover of $W$.

The map $H_3(W;\Lambda)\to H_3(W;\Z_2)$ factors through
$H_3(W;\Lambda)\otimes_\Lambda \Z_2$, more precisely the relation 
is given by a universal coefficient spectral sequence
$$\Tor_p^\Lambda(H_q(W;\Lambda),\Z_2)\Rightarrow H_{p+q}(W;\Z_2).$$
(This can also be interpreted as the Serre spectral sequence for the fibration
$\tilde{W}\to W\to \R \Pro$.)
Here the zeroth and the second row are related by non-trivial differentials:
we compare with the corresponding situation for the space $Q$.
As a result, $H_3(W;\Z_2)\cong H_3(W;\Lambda)\otimes_\Lambda \Z_2$.

By Poincar\'e duality, $H_3(W;\Lambda)\cong H^3(W;\Lambda)$,
and this is free over $\Z$, as there is no $\Z$-torsion in $H_2(W;\Lambda)$.
Since $H_3(W;\Lambda)$ is free over $\Z$, it is a sum of summands of the form
$\Lambda$, $\Z_+$ and $\Z_-$ \cite{CR}.
We also get that $H^3(W;\Lambda)\cong Hom_\Lambda(H_3(W;\Lambda),\Lambda)$, for example again
from a universal coefficient spectral sequence.
The map 
$H_3(W;\Lambda)\to H^3(W;\Lambda) \to Hom_\Lambda(H_3(W;\Lambda),\Lambda)$
describes the $\Lambda$-valued intersection form on $H_3(W;\Lambda)$.
The $\Z_2$-valued intersection form on $H_3(W;\Z_2)$ is given by tensoring with $\Z_2$.

But this implies that for a class $x\in H_3(W;\Lambda)$ with $Tx=\pm x$, 
its image in $H_3(W;\Z_2)$ has intersection 0 with all other elements,
hence it must be 0:
if $Tx=\pm x$, then  
$T\lambda(x,y)=\lambda(Tx,y)=\lambda(\pm x,y)=\pm \lambda(x,y)$, so $\lambda(x,y)$ is a multiple of 
$(1 \pm T)$ and its reduction in $\Z_2$ is zero.

Hence $H_3(W;\Lambda)$ is a free $\Lambda$-module with non-degenerate intersection form.
Since we are in dimension 3, there is a quadratic refinement in 
$\Lambda / \langle x+\bar{x} ,1 \rangle$ which is uniquely determined 
by the intersection form.
Hence we obtain an element in $\tilde{L_6}(\Lambda,w=-)=0$, as $L_6\cong \Z_2$ is given by
the Arf invariant.
Thus it is possible (after stabilization) to do surgery which makes $H_3(W;\Lambda)=0$.

The argument shows that every class of $\Omega_6^A$ 
contains a manifold $W$ with normal 2-type $A$ and $H_3(W;\Lambda)=0$.

\bigskip

For the uniqueness result we take two such manifolds $W_0,W_1$, 
and assume there is a normal $A$-bordism between them.
By Kreck's theory \cite{Kre}, p. 734, there is a surgery obstruction 
in $\tilde{l}_7(\Lambda,w=-1)$ for turning the normal $A$-bordism into an $s$-cobordism.
By surgery below the middle dimension we may assume that the bordism $Y$ is equipped with a 
3-connected map to $A$. Now Kreck defines the surgery obstruction using a construction
of a certain disjoint union $U$ of submanifolds of $Y$ diffeomorphic to $S^3\times D^4$, and defines
the surgery obstruction to be the kernel of 
$H_3(\partial U;\Lambda)\to H_3(Y\setminus int(U), W_0;\Lambda)$. 
He also notes that the orthogonal complement of this kernel is the kernel of
$H_3(\partial U;\Lambda)\to H_3(Y\setminus int(U), W_1;\Lambda)$.
But in our case $H_3(W_i;\Lambda)=0$ so that both these kernels are equal to the kernel of
$H_3(\partial U;\Lambda)\to H_3(Y\setminus int(U); \Lambda)$.

This implies that the surgery obstruction just defined lies in the group
$\tilde{L}_7(\Lambda,w=-1)$ which is trivial as computed by Wall.
Hence we get as a result that $A$-bordant manifolds $W_i$ with normal 2-type $A$ and $H_3(W;\Lambda)=0$
are diffeomorphic (relative to $A$).

The proof for normal 2-type $B$ is obtained by replacing all occurrences of $A$ by $B$.
\qed

\subsection{The transfer}

For the transfer (double cover) map $\Omega_6^A\to \Omega_6^{\tilde{A}}\cong \Omega_6^{Spin}(\C\Pro)$
we compare the Atiyah-Hirzebruch spectral sequences.
We fix the isomorphism 
\begin{eqnarray*}
\Omega_6^{Spin}(\C\Pro) & \cong & \Z^2\\
\left[f:M\to \C\Pro\right] & \mapsto &  \left( \left\langle \frac{p_1(M)f^*x-4f^*x^3}{24},[M]\right\rangle,
\langle f^*x^3,[M]\rangle \right).
\end{eqnarray*}
One computes the homology transfers using the long exact sequences coming from short exact coefficient sequences:
\begin{eqnarray*}
H_6(Q;\Z_-)&\stackrel{\cong}{\to}& H_6(\C\Pro), \\
H_4(Q;\Z_2)&\stackrel{0}{\to}& H_4(\C\Pro;\Z_2), \\
H_2(Q;\Z_-)&\stackrel{\cong}{\to}& H_2(\C\Pro).
\end{eqnarray*}
We see that the transfer gives a map of short exact sequences
$$
\xymatrix{
F_5 \cong \Z_2\oplus\Z 
\ar[d]_{(a,b)\mapsto 2b}
\ar[rr]^{(a,b)\mapsto (a,b,0)}
&&
\Omega_6^A \cong \Z_2\oplus \Z^2
\ar[rr]^{(a,b,c)\mapsto c}
\ar[d]^{tr}
&&
E^{\infty}_{6,0}\cong \Z
\ar[d]^{c\mapsto 2c}
\\
\tilde{F}_5 \cong \Z 
\ar[rr]_{b\mapsto (b,0)}
&&
\Omega_6^{\tilde{A}} \cong \Z^2
\ar[rr]_{ (b,c)\mapsto c}
&&
\tilde{E}^{\infty}_{6,0}\cong \Z
}
$$
which shows (using the snake lemma) that 
$\Omega_6^A\to \Omega_6^{\tilde{A}}$ has a cokernel of order 4.
But the composition of projection and transfer: 
$\Omega_6^{\tilde{A}} \to \Omega_6^A\to \Omega_6^{\tilde{A}}$
is multiplication by 2.
So the image of the transfer consists exactly of all classes divisible by 2.

It also follows that one can find generators for the free summands in $\Omega_6^A$
as images of the generators of $\Omega_6^{\tilde{A}}$.

\bigskip

\begin{thm}\cite{confree}
The image of the transfer map $\Omega_6^B\to \Omega_6^{\tilde{B}}\cong \Omega_6^{Spin}(\C\Pro)$ is $2\Z \oplus \Z$.
\end{thm}

\subsection{Generators for the bordism group $\Omega_6^B$}

We use the Thom isomorphism for twisted spin bordism: 
$\Omega_6^B\cong \Omega_6^{Spin}(Q;L) \cong \Omega_7^{Spin}(DL,SL)$.
Under this isomorphism, the boundary map $\Omega_7^{Spin}(DL,SL)\to \Omega_6^{Spin}$
corresponds to the transfer map $\Omega_6^B\to \Omega_6^{\tilde{B}}$.
It follows that the torsion elements in $\Omega_6^B\cong \Omega_7^{Spin}(DL,SL)$
come from $\Omega_7^{Spin}(DL)\cong \Omega_7^{Spin}(Q)$. Moreover, $\Omega_7^{Spin}(\C \Pro)=0$,
as one sees easily from the Atiyah-Hirzebruch spectral sequence.
Thus $\Omega_7^{Spin}(Q)\cong \Z_4$ must be responsible for the torsion.
Moreover, a study of the Atiyah-Hirzebruch spectral sequence shows that the map 
$\Omega_7^{Spin}(\C \Po^3/\tau)\to \Omega_7^{Spin}(Q)$ is an isomorphism.
Now there is a bundle $\R \Po^2 \to \C \Po^3/\tau \to S^4$. So $\C \Po^3/\tau - \R \Po^2$ is
an $\R \Po^2$-bundle over $D^4$, so homotopy equivalent to $\R \Po^2$.
Again the Atiyah-Hirzebruch spectral sequence shows that 
$\Omega_7^{Spin}(\C \Po^3/\tau - \R \Po^2)=\Omega_6^{Spin}(\C \Po^3/\tau -\R \Po^2)=0$.
Thus $\Omega_7^{Spin}(\C \Po^3/\tau)\cong \Omega_7^{Spin}(\C \Po^3/\tau, \C \Po^3/\tau - \R \Po^2)$,
and the latter is isomorphic to $\Omega_3^{Spin}(\R \Po^2)\cong \Z_4$ by the Thom isomorphism.
Again the Atiyah-Hirzebruch spectral sequence shows that the transfer 
$\Omega_3^{Spin}(\R \Po^2)\to \Omega_3^{Spin}(S^2)\cong \Z_2$ is surjective.
Thus one can detect a generator of $\Omega_7^{Spin}(\C \Po^3/\tau)$ by the composition
$$\Omega_7^{Spin}(\C \Po^3/\tau)\to \Omega_7^{Spin}(\C \Po^3/\tau, \C \Po^3/\tau - \R \Po^2)
\to \Omega_3^{Spin}(\R \Po^2)\to \Omega_3^{Spin}(S^2).$$

Now take the seven-dimensional manifold $(\C \Po^3 \times S^1)/(\tau,c)$, where $c$ is complex conjugation.
The spin structure on $\C \Po^3 \times S^1$ which restricts to the non-bounding one on $S^1$ is
preserved by the involution $(\tau,c)$, so that we obtain a spin structure on the quotient.
The map to $\C \Po^3/\tau$ is projection on the first coordinate.
This intersects $\R \Po^2=\C \Po^1/\tau$ transversely, so that the Thom isomorphism maps it to the pullback
$(\C \Po^1\times S^1)/(\tau,c)\to \C \Po^1/\tau$ whose double cover is the projection
$\C \Po^1 \times S^1 \to \C \Po^1$, which is a generator for $\Omega_3^{Spin}(S^2)$ since the Spin structure
restricts to the non-boundant one on $S^1$.

We have to apply the Thom isomorphism $\Omega_7^{Spin}(DL|_{\C \Po^3/\tau},SL|_{\C \Po^3/\tau})\to \Omega_6^{Spin}(\C \Po^3/\tau;L)$ to the
map $(\C \Po^3 \times S^1)/(\tau,c) \to \C \Po^3/\tau \to (\C \Po^3 \times D^1)/(\tau,-1)$. 
For this, we homotope the map to make it transversal to the zero section: take
\begin{eqnarray*}
(\C \Po^3 \times S^1)/(\tau,c) & \to & (\C \Po^3 \times D^1)/(\tau,-1) \\
 \ [ x , y ] & \mapsto & [ x,Im(y) ]
\end{eqnarray*}
and intersect it with the zero section: 
we obtain two copies of $\C \Po^3/\tau$ with different $B$-structures.

It follows that the torsion $\Z_4$ in $\Omega_6^B$ is generated by the sum (or the difference) of two
copies of $\C \Po^3/\tau$ with different $B$-structures. 
Let us denote one of them by $\sigma$ and the other by $\sigma'$.

Thus we obtain as generators for $\Omega_6^B\cong \Omega_7^{Spin}(DL,SL)$:
\begin{itemize}
\item 
$X\times D^1 \stackrel{f\times id}\to \C \Pro \times D^1$, where $X$ is a simply-connected spin 6-manifold with $H^3(X)=0$, 
trivial cup product $H^2(X)\times H^2(X)\to H^4(X)$, and $f:X\to \C \Pro$ defines a generator of $H^2(X)$ such that the first
Pontrjagin class of $X$ is equal to 24 times the generator of $H^4(X)$ which is dual to $f\in H^2(X)$. 
\item
$((\C \Po^3\times D^1)/(\tau,-1),\sigma)$,
\item
$((\C \Po^3\times D^1)/(\tau,-1),\sigma) -((\C \Po^3\times D^1)/(\tau,-1),\sigma')$.
\end{itemize}
The two former generators generate free summands, the latter generates the torsion summand $\Z_4$.
In this basis, the map $\Omega_6^B \cong \Z^2 \oplus \Z_4 \to \Z^2\cong \Omega_6^{Spin}(\C\Pro)$ is given as $(a,b,c)\mapsto (2a,b)$.

\subsection{The automorphism groups of $A$ and $B$ and their action on the bordism groups}

The automorphism group $Aut(B)_{BO}$ of fiber homotopy classes of fiber homotopy self-equivalen\-ces of $B$ acts
on $\Omega_6^B$.
We saw that $\Omega_6^B\cong \Omega_7^{Spin}(DL,SL)$, where 
$L$ is the nontrivial real line bundle $(\C \Pro \times \R)/(\tau,-1)$.

\begin{lem}
The set of equivariant oriented diffeomorphism classes of free involutions on 6-manifolds with $H_3=0$ and
whose quotient spaces have normal 2-type $B$ are given as the orbits.
Again, $B$ can be replaced by $A$ in the theorem.
\end{lem}
Proof: Equivariant diffeomorphism classes of free involutions are the same as diffeomorphism classes of the
quotients. Now the theorem follows from the uniqueness of the Postnikov decomposition, i.e. for a given
manifold $W$ with normal 2-type $B$, the map $W\to B$ is uniquely determined up to fiber homotopy self equivalences of
$B$ over $BO$. See also \cite{KrM}.
\qed
\bigskip

The restriction of the first component of a fiber homotopy self-equivalence 
of $Q\times BSpin$ to $Q$ is a self-homotopy equivalence of $Q$. There is a unique
free homotopy class of maps $Q\times BSpin\to Q$ which is an isomorphism on $\pi_1$,
and thus (using obstruction theory to extend homotopies from $Q$ to $Q\times BSpin$) a unique free 
homotopy class of maps $Q\times BSpin\to Q$ which are an isomorphism on $\pi_1$.
The vertical homotopy classes of fiber homotopy equivalences thus correspond  
to the different choices (up to homotopy) of a homotopy $\eta$ from 
$Q\times BSpin \to Q \to K(\Z_2,1)\times K(\Z_2,2)$
to $Q\times BSpin\to BO \to K(\Z_2,1)\times K(\Z_2,2)$.
So the group has four elements, and the action of the group on the set of normal $B$-structures 
on a manifold $f:M\to B$ (i.e. spin structures on $\nu_M-f^*L$) is given by 
just changing the spin structure $\sigma$  
into $\sigma, -\sigma, \sigma +f^*t, -\sigma+f^*t$, where $t\in H^1(B;\Z_2)$ is the generator.

On $\Omega_7^{Spin}(DL,SL)$, the group $Aut(B)_{BO}$ acts in the following way:
the negative spin generator acts by -1, and the spin flip generator acts by 
$$((\C \Po^3\times D^1)/(\tau,-1),\sigma)\mapsto ((\C \Po^3\times D^1)/(\tau,-1),\sigma'),$$ 
and is the identity on $X$.
Thus in the decomposition $\Omega_7^{Spin}(DL,SL)\cong \Z\oplus \Z \oplus \Z_4$ given by the above
generators, the negative spin generator acts by $(a,b,c)\mapsto (-a,-b,-c)$, 
and the spin flip generator acts by $(a,b,c)\mapsto (a,b,b-c)$.
We obtain orbits of the form $$\{ (a,b,c), (a,b,b-c),(-a,-b,-c), (-a,-b,-b+c)\}.$$

The group $Aut(A)_{BO}$ also has four elements which act 
on the set of normal $A$-structures 
on a manifold $f:M\to A$ (i.e. spin structures on $\nu_M-f^*(3L)$) by 
just changing the spin structure $\sigma$  
into $\sigma, -\sigma, \sigma +f^*t, -\sigma+f^*t$, where $t\in H^1(A;\Z_2)$ is the generator.

This is either the identity or minus the identity on the generators for the free summands in $\Omega_6^A$
as they are images of the generators of $\Omega_6^{\tilde{A}}$ unde projection.
All group elements must act by the identity on the torsion generator.
We obtain orbits $$\{ (a,b,c), (-a,-b,-c)\}.$$

The group $Aut(K(\Z,2) \times BSpin)_{BO}$ also has four elements,
which act on the bordism group by reversing the spin structure and/or the class in $H^2$.
It follows that here the orbits are of the form $$\{ (a,b), (-a,-b) \}.$$

\bigskip

For the proof of theorem \ref{Mthm}, it is now sufficient to count preimages and orbits:
An element of the form $(2a+1,b)\in \Omega_6^{Spin}(\C \Pro)$ has no preimages in $\Omega_6^A$ or $\Omega_6^B$.
An element $(2a,2b+1)$ in $\Omega_6^{Spin}(\C \Pro)$ has four preimages $(a,2b+1,c)$ in $\Omega_6^B$
and no preimages in $\Omega^A_6$. These four preimages, together with the four preimages of $(-2a,-2b-1)$,
decompose into two orbits.
For $(2a,2b)$ in $\Omega_6^{Spin}(\C \Pro)$ we obtain four preimages $(a,2b,c)$ in $\Omega_6^B$.
These four preimages, together with the four preimages of $(-2a,-2b)$, decompose into three orbits.
The element $(2a,2b)$ in $\Omega_6^{Spin}(\C \Pro)$ has two preimages $(a,b,c)$ in $\Omega_6^B$.
These two preimages, together with the two preimages of $(-2a,-2b)$, decompose into two orbits.

\section{Involutions on $S^6$ with three-dimensional fixed point set} 

\subsection{Non-equivariant classification of embeddings in $S^6$}

Before we classify involutions on $S^6$, let us recall the non-equivariant results on embeddings
of three-manifolds into $S^6$.

Naturally the most interesting case is the one of knotted 3-spheres in the six-sphere.
Here the results are due to Haefliger.
One could consider various equivalence relations on knotted $S^3$'s in $S^6$.
In all cases there is an addition defined using connected sums: 

First we can look at the group $C_3^3$ of isotopy classes of embeddings of $S^3$ into $S^6$.
Second, the group $\Theta$ of diffeomorphism classes of embeddings of $S^3$ into oriented manifolds $X$ diffeomorphic to 
$S^6$, relative to $S^3$. (We require a diffeomorphism to be the identity on $S^3$ and to preserve orientations.)
Third, the group $\Sigma$ of orientation-preserving diffeomorphism classes of pairs $(X,M)$, 
where the oriented manifold $X$ is diffeomorphic to $S^6$ 
and the oriented submanifold $M$ is diffeomorphic to $S^3$.

There are obvious surjective group homomorphisms $C_3^3\to \Theta\to\Sigma$.
Haefliger showed that $C_3^3$ and $\Sigma$ are both isomorphic to $\Z$, so
that both of the above maps are isomorphisms. 

To prove that $C_3^3\to \Theta$ is an isomorphism, one needs to show that a diffeomorphism $S^6\to S^6$ relative the embedded $S^3$
can be replaced by an isotopy. This is true since $\pi_0(Diff(D^n,\partial))\to \pi_0(Diff(S^n))$ is surjective. This means 
it is always possible to modify the original diffeomorphism on a disk such that the resulting diffeomorphism is isotopic to the identity.

One explanation of the isomorphism $\Theta\to\Sigma$ is Cerf's result that $Diff^+(S^3)$ is connected. 
Thus every orientation-preserving diffeomorphism of $S^3$ is isotopic to the identity,
and this isotopy extends to an ambient equivariant isotopy of $S^6$. 

The negative of an isotopy class is given by precomposing the embedding $S^3\to S^6$ with an 
orientation-reversing self-diffeomorphism of $S^3$.

\bigskip

Isotopy classes of framed embeddings $S^3\times D^3\to S^6$ are in bijection with $\Z^2$, the framing giving an additional integer invariant
(the isomorphism to $\Z^2$ depends on a choice). 
Note also that in the PL category, all non-framed knots are trivial, but the isotopy classes of framed knots are in bijection
with smooth isotopy classes of smooth framed knots.

\bigskip

For embeddings of general closed oriented connected 3-manifolds $M^3$ into $S^6$,
the argument that diffeomorphism relative to the submanifold implies isotopy of the embeddings still holds. 
Isotopy classes of embeddings $Emb(M,S^6)$ have been classified by Skopenkov in \cite{Sko}.
To an isotopy class of embeddings $i:M\to S^6$ he associates its Whitney invariant $W(i)\in H_1(N)$. (For the precise definition
we refer to \cite{Sko}, we give a description in special cases in remark \ref{WhitInv}.)
The map $W:Emb(M,S^6)\to H_1(N)$ is surjective, and $C^3_3\cong \Z$ acts transitively 
on the fibers by connected sum. This action has non-trivial stabilizers in general: 
There is a bijective map (the Kreck invariant) from the fiber over $u\in H_1(N)$ 
to $\Z_{d(u)}$, where $d(u)$ is the divisibility of $\bar{u}\in H_1(M)/\{\text{torsion}\}$. 
In general, both the Whitney and the Kreck invariant depend on choices.
For $\Z_2$-homology spheres $M$, the Whitney invariant does not depend on choices, and $C_3^3$ acts freely on the fibers,
so that $Emb(M,S^6)$ is in non-canonical bijection with $\Z \times H_1(M)$. 
Instead of a map to $\Z$, the Kreck invariant describes an action of $\Z$ on $Emb(M,S^6)$ which leaves the Whitney invariant fixed.

\subsection{What should we classify in the equivariant case?}

In the equivariant case, there are again various equivalence relations one can put on the set of embeddings respectively involutions.
However, in order to get a well-defined connected sum operation, it is necessary to orient both the 6-manifold and the fixed point set.

\begin{remark}
The proof of the uniqueness of the non-equivariant connected sum construction 
can be generalized to show that the equivariant diffeomorphism type of the 
equivariant connected sum of two conjugation manifolds of dimension $2n$ 
only depends on the connected component of the set of equivariant isomorphisms of a 
tangent space at a fixed point with $(\R^{2n},(1^n,-1^n))$. More generally, varying the chosen fixed point,
we get a bundle of such isomorphisms over the fixed point set, and the equivariant diffeomorphism type 
depends only on the connected component in the total space of this bundle. 
(The total space has two components if the fixed point set is orientable, and one component if it is not.)
See also Definition 1.1 and Lemma 1.2 of \cite{Loe}.
In particular, the connected sum is unique up to equivariant diffeomorphism
if we are provided with orientations of the conjugation manifolds and their fixed point sets,
but in general depends also on an orientation of the fixed point sets.
This answers a question in \cite{HHP}, and it also explains why we are less interested in
the classification of conjugations without an orientation of the fixed point set. 
\end{remark}

We fix $M^3$, allow various involutions $\tau$ on $S^6$, and
consider equivariant embeddings $(M,id)\to (S^6,\tau)$ such that the image of $i$ is the fixed point set of $\tau$.
Again there are several equivalence relations which we can put on this set.
\begin{itemize}
\item
Two embeddings $i_0:(M,id)\to (S^6,\tau_0)$ and $i_1:(M,id)\to (S^6,\tau_1)$
are equivalent if there is an equivariant, orientation-preserving diffeomorphism  $\phi:(S^6,\tau_0)\to (S^6,\tau_1)$ relative $M$,
i.e. there is a commutative triangle 
$$
\xymatrix{
(M,id) 
\ar[r]^{i_0}
\ar[dr]^{i_1}
& (S^6,\tau_0)
\ar[d]^{\phi}
\\
& (S^6,\tau_1).
}
$$
This classsifies involutions together with an identification
of the fixed point set with the given 3-manifold $M$.
We get a set $Emb_{\Z_2}(M,S^6)$.
\item
Two embeddings $i_0:(M,id)\to (S^6,\tau_0)$ and $i_1:(M,id)\to (S^6,\tau_1)$
are equivalent if there is an equivariant diffeomorphism  $\phi:(S^6,\tau_0)\to (S^6,\tau_1)$ which restricts to some self-diffeomorphism
$\phi_M$ of $M$. We require that both $\phi$ and $\phi_M$ are orientation-preserving. There is a commutative square 
$$
\xymatrix{
(M,id) 
\ar[r]^{i_0}
\ar[d]^{\phi_M}
& (S^6,\tau_0)
\ar[d]^{\phi}
\\
(M,id) 
\ar[r]^{i_1}
& (S^6,\tau_1).
}
$$
This classifies (up to orientation-preserving diffeomorphism) involutions plus an orientation of the fixed point set.
One might call this the classification of conjugations up to conjugation.
We get a set $Inv_M(S^6)$. 
\end{itemize}

\begin{remark}
One could also define an equivariant version of isotopy classes of embeddings:
we say that two equivariant embeddings $i_0:(M,id)\to (S^6,\tau_0)$ and $i_1:(M,id)\to (S^6,\tau_1)$
are equivalent if there is an equivariant embedding 
$i:(M\times I,id)\to (S^6\times I,\tau)$ 
such that $i(x,t)=(i_t(x),t)$ and $\tau(y,t)=(\tau_t(y),t)$, 
hence in particular it restricts on both ends to $i$ respectively $i'$.
(We also require that the image of all the embeddings involved is the whole fixed point set.)
\end{remark}

\subsection{Analysis of involutions on $S^6$ with three-dimensional fixed point set} 

We recall the classical result:
\begin{thm}[P.A.~Smith \cite{Bre}]
If an involution on $S^n$ has fixed points, the fixed point set is a $\Z_2$-homology sphere.
\end{thm} 
Comparing with Skopenkov's classification of embeddings, our first result is the following.
\begin{pro}
If $i:M^3\to S^6$ is the embedding of a fixed point set of an involution, then the Whitney invariant of the embedding vanishes: $W(i)=0$.
\end{pro}
Proof: Let $\tau$ be such an involution, and let $\sigma$ be the reflection of $S^6$ at the equator.
In \cite{Sko} it is proved that $W(\sigma \circ i)=-W(i)$. 
But we have a commutative square
$$ 
\xymatrix{
M 
\ar[d]^=
\ar[r]^{i}
& S^6
\ar[d]^{\sigma \circ \tau}
\\
M
\ar[r]^{\sigma \circ i}
& S^6
}
$$
which shows that $i$ and $\sigma \circ i$ are isotopic. Thus $W(i)=-W(i)$, and since $H_1(M)$ consists of odd torsion 
only, we have $W(i)=0$. \qed

\begin{remark}\label{WhitInv}
Up to a multiplicative factor, which is invertible in the case of $\Z_2$-homology spheres, 
the Whitney invariant can also be defined in the following way:
Let $C_i=S^6\setminus i(M^3)$. By Alexander duality, $H^2(C_i)\cong \Z$ and $H^4(C_i)\cong H^2(M)\cong H_1(M)$.
Then the invariant is given by the square of a generator of $H^2(C_i)$.
In the case of a fixed point set of a conjugation, the involution acts by multiplication with $-1$ on both $H^2(C_i)$ and $H^4(C_i)$.
It follows that the square of a generator of $H^2(C_i)$ must be 0.
\end{remark}

\bigskip

Recall that by the equivariant tubular neigborhood theorem, for every involution on $S^6$ with fixed point set $M$ 
we can write $S^6=M\times D^3 \cup V$, such that the involution is $id\times -id$ on $M\times D^3$, and free on $V$.
The quotient $W=V/\tau$ is a manifold with fundamental group $\Z_2$ and boundary $M\times \R \Po^2$.
The normal 2-type of $W$ is $B$.
The inclusion of the boundary needs to induce an isomorphism $\pi_2(\R \Po^2)\to \pi_2(B)$.

In order to define the bordism set $\Omega_6^{(B,M\times \R \Po^2)}$, 
we have to fix a normal $B$-structure on $M\times \R \Po^2$. Since we are interested in a classification,
we first consider all relevant normal $B$-structures.

For any such structure, the homotopy class of maps $M\times \R \Po^2\to Q\to \R \Pro$ 
is the non-trivial class in $H^1(M\times \R \Po^2;\Z_2)$. We have to lift this non-trivial map 
$M\times \R \Po^2\to \R \Pro$ to $Q$. Any lift, together with a choice 
of a spin structure on $\nu_{M\times \R \Po^2} - f^*(L)$, will then be a normal $B$-structure on $M\times \R \Po^2$. 
It is easy to find a lift $f$ on $\R \Po^2$, see \cite{OlbThesis}, p.47, and composing any lift with the projection
$M\times \R \Po^2\to \R \Po^2$ gives a lift $S^3\times \R \Po^2\to Q$. 
Obstruction theory shows that pointed homotopy classes of lifts on $\R \Po^2$ are classified by 
$H^2(\R \Po^2,\Z_-)\cong \Z$, and that every lift on $\R \Po^2$ extends uniquely up to homotopy
to a lift on $M\times \R \Po^2$, see \cite{OlbThesis}, pp.50-52.

But there is a further condition. Only two pointed homotopy classes of lifts induce an 
isomorphism $\pi_2(\R \Po^2)\to \pi_2(Q)$,
and one obtains one from the other by precomposing with the nontrivial pointed homotopy class
of maps $\R \Po^2 \to \R \Po^2$.
However, this map is freely homotopic to the identity of $\R \Po^2$.
Thus we get up to (free) homotopy a unique map $f:\R \Po^2\to Q$.
There are four spin structures on $\nu_{M\times \R \Po^2} - f^*(L)$, 
since spin structures on unoriented (but orientable and spin) bundles over $N$ 
are in bijection with $H^0(N;\Z_2)\times H^1(N;\Z_2)$.
Thus there are four distinguished normal $B$-structures on $M \times \R \Po^2$ which can be used in the construction. 

\begin{thm}
For every distinguished normal $B$-structure on $M\times \R \Po^2$, 
every element of the bordism set $\Omega_6^{(B,M\times \R \Po^2)}$ contains (up to diffeomorphism relative to $B$ and the boundary) a unique manifold $W$
which produces a conjugation on $S^6$.
\end{thm}
\begin{remark}
Actually we also see that $H_3(W;\Lambda)$ consisting just of odd torsion is a necessary and sufficient condition for this.
\end{remark}
Proof: The proof of existence is a slight modification of the proof of theorem \ref{freebord}.
(For full details, see the proof of theorem 1.3 in \cite{confree}.)
Also the uniqueness extends from the proof of theorem \ref{freebord}:
If we take two manifolds $W_0,W_1$ with the same normal $B$-structure on the boundary,
and which both produce conjugations on $S^6$, and assume that there is a normal $B$-bordism between them,
we have to modify the argument from \ref{Con} slightly to show that Kreck's surgery obstruction is 0:
For both $i=0,1$ the kernel of $H_3(\partial U;\Lambda)\to H_3(Y\setminus int(U),W_i;\Lambda)$ is equal to the kernel
of $H_3(\partial U;\Lambda)\to H_3(Y\setminus int(U);\Lambda) /\{\text{torsion}\}$. So again the surgery obstruction
lies in $\tilde{L}_7(\Lambda,w=-1)=0$.
\qed

\begin{cor}
Diffeomorphism classes of $W$ relative to $M\times \R \Po^2$ and to $B$ 
(where $M\times \R \Po^2 \to B$ is fixed) producing conjugations on $S^6$
are in bijection with $\Omega_6^{(B,M\times \R \Po^2)}\cong \Z^2\oplus \Z_4$.
\end{cor}

Diffeomorphism classes of $W$ relative to $M\times \R \Po^2$ and to $B$ 
(where $M\times \R \Po^2 \to B$ is not fixed) producing conjugations on $S^6$
are in bijection with the disjoint union of these four bordism sets.

\bigskip

Equivariant connected sum of two conjugations with fixed point sets $M_1$ respectively $M_2$ corresponds
to a map of bordism sets defined by gluing along part of the boundary (or parametrized boundary connected sum):
Choose disks $D^3$ in $M_1,M_2$ centered at the points where the connected sum is performed.
Then there is a map
\begin{eqnarray*}
\Omega_6^{(B,M_1\times \R \Po^2)} \times \Omega_6^{(B,M_2\times \R \Po^2)} &\to& \Omega_6^{(B,(M_1\# M_2)\times \R \Po^2)} \\
 (W_1, W_2) & \mapsto & W_1\cup_{D^3\times \R \Po^2} W_2
\end{eqnarray*}
This is equivariant with respect to the action of the bordism group $\Omega_6^B$.
Hence it equips $\Omega_6^{(B,S^3\times \R \Po^2)}$ with a group structure, and we obtain an action of 
$\Omega_6^{(B,S^3\times \R \Po^2)}$ on $\Omega_6^{(B,M\times \R \Po^2)}$ for any $M$.

\bigskip

To compare with the non-equivariant embedding results of Skopenkov, it suffices to consider the image under 
the transfer map:  
Comparing with \cite{Sko}, we see that for the embeddings with Whitney invariant 0, 
the Kreck invariant describes an action of $\Z$ on $Emb(M,S^6)$
which corresponds precisely to the action of 
$$\Omega_6^{Spin}(K(\Z,2),\partial=S^3\times S^2)/(0\oplus\Z\subset \Omega_6^{Spin}(K(\Z,2)))$$ on
the set $$\Omega_6^{Spin}(K(\Z,2),\partial=M\times S^2)/(0\oplus\Z\subset \Omega_6^{Spin}(K(\Z,2))).$$ 
We will see in the next section how the quotients arise also in the equivariant setting. 

\subsection{Equivariant diffeomorphism classes as a quotient by group actions}

Now we have to relate the set of relative diffeomorphism classes of the previous section to 
the set $Emb_{\Z_2}(M,S^6)$ of equivariant diffeomorphism classes of six-spheres with involution whose 
fixed point set is identified with $M$.

Basically, we forget the $B$-structure, we construct 
the equivariant inclusion $M\times D^3\to X$ from $M\times \R \Po^2 \to W$,
and we forget the tubular neighbourhood and the framing of the normal bundle.

More precisely, we have the following sets of equivalence classes:
\begin{enumerate}
\item \label{item1}
The set $T_1$ of diffeomorphism classes of manifolds $W$ relative to $M\times \R \Po^2$ and $BO$.
An element is represented by $M\times \R \Po^2 \to W \to B$, such that the first map is the inclusion of the boundary. 
Two representatives $W$ and $W'$ are equivalent if
there is a diffeomorphism $W\to W'$ which commutes with the maps from $M\times \R \Po^2$ and to $BO$.
\item \label{item2}
The set $T_2$ of diffeomorphism classes of manifolds $W$ relative to $M\times \R \Po^2$.
An element is represented by $M\times \R \Po^2 \to W$, which is the inclusion of the boundary. Two representatives $W$ and $W'$ are equivalent if
there is a diffeomorphism $W\to W'$ which commutes with the maps from $M\times \R \Po^2$.
\item
The set $T_3$ of equivariant diffeomorphism classes of manifolds $V$ relative to $M\times S^2$.
An element is represented by $(M\times S^2,(id,-id)) \to (V,\tau)$, which is the inclusion of the boundary, and where $\tau$ is a free involution.
Two representatives $V$ and $V'$ are equivalent if
there is a diffeomorphism $(V,\tau)\to (V',\tau')$ which commutes with the maps from $M\times S^2$.
\item
The set $T_4$ of equivariant diffeomorphism classes of closed manifolds $X$ relative to $M\times D^3$.
An element is represented by an equivariant embedding $(M\times D^3,(id,-id)) \to (X,\tau)$, where $\tau$ is free on the complement of the image.
Two representatives $X$ and $X'$ are equivalent if
there is a diffeomorphism $(X,\tau)\to (X',\tau')$ which commutes with the maps from $M\times D^3$.
\item
The set $Emb_{\Z_2}(M,S^6)$ of equivariant diffeomorphism classes of manifolds $X$ relative to $M$.
An element is represented by an equivariant embedding $(M,id) \to (X,\tau)$, where $\tau$ is free on the complement of the image. 
Two representatives $X$ and $X'$ are equivalent if
there is a diffeomorphism $(X,\tau)\to (X',\tau')$ which commutes with the maps from $M$.
\end{enumerate}

On the set $T_1$, the automorphism group $Aut(B)_{BO}$ of fiber homotopy classes of fiber homotopy 
self-equivalences of $B$ acts by post-composing $M\times \R \Po^2 \to B$ with $B\to B$.

We saw that the group $Aut(B)_{BO}$ has four elements, and the action of the group on the set of normal $B$-structures 
on a manifold $f:M\to B$ (i.e. spin structures on $\nu_M-f^*L$) is given by 
just changing the spin structure $\sigma$  
into $\sigma, -\sigma, \sigma +f^*t, -\sigma+f^*t$, where $t\in H^1(B;\Z_2)$ is the generator.

Thus the group $Aut(B)_{BO}$ acts freely and transitively on the set of 
distinguished normal $B$-structures on $M\times \R \Po^2$. In particular
each orbit of the action on the set in \ref{item1} contains a unique element from each of the four bordism sets.

By uniqueness of the Postnikov decomposition, forgetting the map to $B$ is a bijection from 
the quotient of $T_1$ by $Aut(B)_{BO}$ to the set $T_2$. See also \cite{KrM}.
Thus the set $T_2$ is in bijection with $\Omega_6^{(B,M\times \R \Po^2)}\cong \Z^2\oplus \Z_4$.

\bigskip

The sets $T_2$, $T_3$ and $T_4$ are in bijective correspondence:
One gets from $M\times D^3 \to X$ by restriction to $M\times S^2 \to V=X\setminus int(M\times D^3)$
and from $M\times S^2 \to V$ one takes the quotient by the involution to obtain $M\times \R \Po^2 \to W=V/\tau$.
Vice-versa, from $M\times \R \Po^2 \to W$, take any non-trivial double covering $V$ of $W$ and any map
$M\times S^2 \to V$ inducing the given $M\times \R \Po^2 \to W$. The involution on $V$ is the non-trivial deck transformation.
Up to equivalence, this does not depend on the choices made. 
And from $\phi:M\times S^2\to V$, obtain $M\times D^3\to X=M\times D^3\cup_\phi V$. The smooth structure on the latter 
depends on choices, but only up to equivalence.

\bigskip

On the sets $T_2$, $T_3$ and $T_4$, we have an action of the group of bundle automorphisms $Map(M,O(3))$ of $M\times D^3$: 
every bundle automorphism induces a self-diffeomorphism of the boundary $M\times S^2$ 
which induces a self-diffeomorphism of $M\times \R \Po^2$, and we precompose with 
these diffeomorphisms. 
A bundle automorphism $f:M\to O(3)$ corresponds to the diffeomorphism 
\begin{gather*}
\phi_f:M\times \R \Po^2 \to M\times \R \Po^2\\
(x,\pm y)\mapsto (x,\pm f(x)\cdot y)
\end{gather*}
for $y\in S^2$.
So the bundle automorphism which is minus the identity on each fiber acts trivially on
$M\times \R \Po^2$, hence the action is trivial also on the other sets.  
(This corresponds in $T_4$ to the fact that the embeddings $i:M\times D^3 \to X$ and
$M\times D^3\stackrel{id\times -id}\to M\times D^3 \stackrel i\to X$ are related by the equivariant diffeomorphism 
$\tau:(X,\tau)\to (X,\tau)$.)
Thus we may restrict to orientation preserving bundle automorphisms.

\bigskip

The equivariant tubular neighbourhood of the fixed point set is unique up to isotopy
and bundle automorphisms. But isotopies of embeddings can be enlarged to isotopies 
of the ambient space, so that they act trivially on the set of equivariant diffeomorphism classes.
It follows that the action of the group of bundle automorphisms descends to an action of $[M,SO(3)]$, and
that dividing out this action, we get the set $Emb_{\Z_2}(M,S^6)$.

\bigskip

The part which is a little more complicated is to see the action of $\Z \cong \pi_3(SO(3))\cong [M,SO(3)] \ni f$ on a bordism set
$\Omega_6^{(B,M\times \R \Po^2)}$.
We can also precompose with the corresponding self-diffeomorphism $\phi_f$ of $M\times \R \Po^2$, but this changes the precise map
$M\times \R \Po^2 \to B$. Still the normal $B$-structure is preserved:
Up to homotopy, we may assume that $f:M\to SO(3)$ is trivial on a whole disk $D^3$, so that the normal $B$-structure on
$D^3\times \R \Po^2$ does not change. Spin structures on vector bundles over $M\times \R \Po^2$
are in bijection with spin structures on their restrictions to $\R \Po^2$ (using the natural framing 
of the normal bundle), and the latter are invariant. 
Hence these bundle automorphisms preserves the normal $B$-structure on $M\times \R \Po^2$, and we may consider the action of $f$
on the bordism group as given by gluing a mapping cylinder of $f$.

\begin{lem}
Let $B\to BO$ be a fibration, let $N$ be an $(n-1)$-manifold with normal $B$-structure,
$W_1, W_2, W_3$ be normal $B$-nullbordisms of $N$, let $\phi:N\to N$ be a diffeomorphism preserving the normal $B$-structure, let
$C_\phi$ be the mapping cylinder of $\phi$, let $i_0,i_1:N\to C_{\phi}$ be the two natural inclusions, and let $T_\phi$ be the mapping torus 
of $\phi$.
Then, in the bordism group $\Omega_n^B$, we have 
$$(W_1 \cup_{i_0}  C_\phi \cup_{i_1} - W_2) - (W_1 \cup_{id_N} - W_2) = T_{\phi} = W_3 \cup_{i_0}  C_\phi \cup_{i_1} -W_3.$$
\end{lem}
Proof: $W\times I$ can be considered as a bordism between the manifolds
with boundary $\partial W \times I$ and $-W \cup W$. 
The second equality in the statement follows from applying this to $W=W_3$.
Similarly, the first equality is obtained by applying this to $W=W_1$ and $W=W_2$. \qed

\bigskip

It follows that the action of $f$ on the bordism set $\Omega_6^{(B,S^3\times \R \Po^2)}$ is the same as 
taking the disjoint sum with $D^4\times \R \Po^2 \cup_{\phi_f} D^4\times \R \Po^2$. 
The latter is an $\R \Po^2$-bundle on $S^4$, and the corresponding double cover is a $S^2$-bundle 
over $S^4$ which can be identified with $S^2\to \C \Po^3 \to \HB \Po^1$.
Thus the induced action on $\Omega_6^{(\tilde{B},S^3\times S^2)}\cong\Z^2$ is by a generator
for the first summand of $\Omega_6^{\tilde{B}}\cong \Z^2$.

\bigskip

It follows that the set $Emb_{\Z_2}(M,S^6)$ of orbits of the action on $T_4$ is in bijection with $\Z \oplus \Z_4$.
In particular $Emb_{\Z_2}(M,S^6)$ is a group which acts freely and transitively on $Emb_{\Z_2}(M,S^6)$ for all $M$.

Comparing with the non-equivariant classification of embeddings in $S^6$, we see that forgetting the involution
defines a map $Emb_{\Z_2}(M,S^6) \cong \Z \oplus \Z_4 \to Emb(M,S^6) \cong \Z\oplus H_1(M)$ which is equivariant with respect to the 
group homomorphism $Emb_{\Z_2}(S^3,S^6) \cong \Z \oplus \Z_4 \to Emb(S^3,S^6) \cong \Z$ given by $(a,b)\mapsto 2a$.
In particular the image of $Emb_{\Z_2}(M,S^6)$ is acted upon freely by $2\Z\subseteq \Z$, and the map is 4-to-1.

For embeddings of $\Z_2$-homology spheres we saw that the elements $[i:M\to S^6] \in Emb(M,S^6)$ with vanishing Whitney invariant
are acted freely and transitively upon by $C_3^3\cong \Z$. The subset of isotopy classes of embeddings which are the fixed point 
sets of conjugations are acted freely and transitively upon by $2\Z\subseteq \Z$. There are up to equivariant diffeomorphism
relative to $i$ exactly four such conjugations for every $i$. This proves theorems \ref{thm13} and \ref{embcl}.

\subsection{From embeddings to submanifolds - proof of theorem \ref{invcl}}
The more natural thing is to classify involutions without the additional identification of the fixed point set with 
a fixed 3-manifold $M$. The invariant of the involution should be its fixed point set, i.e. a submanifold of $S^6$.
In order to get from embeddings up to diffeomorphism to submanifolds up to diffeomorphism, it suffices to divide out the action
of the group of self-diffeomorphisms $Diff(M)$.
Since isotopies extend to ambient isotopies (this also holds in this equivariant case, since it suffices to extend a vector field on the
fixed point set to an equivariant vector field on the whole space), 
and these give equivariant diffeomorphisms, the action of $Diff(M)$ factors through the
mapping class group of $M$.

Since in our case, the map $M\times \R \Po^2 \to Q$ factors through $\R \Po^2$, this map does not change, so that we get
the same normal $B$-structure. Then the action of a self-diffeomorphism $f:M\to M$ on the bordism set
$\Omega_6^(B,M\times \R \Po^2)$ is by disjoint union with the mapping torus $T_{f\times id}$ of $f\times id: M\times \R \Po^2 \to \R \Po^2$.
Now $T_{f\times id}=T_f \times \R \Po^2$, and the map to $Q$ factors again through $\R \Po^2$. Spin structures twisted by $L$ on
$T_f\times \R \Po^2$ are products of a spin structure on $M$ and a spin-structure twisted by $L$ on $\R \Po^2$. Thus a spin-nullbordism 
of $T_f$ gives a normal $B$-nullbordism of $T_{f\times id}$. Recall that a four-dimensional spin manifold is zero bordant iff its 
signature is zero, and that the signature of a mapping torus is always zero. Hence the action of the mapping class group of $M$ on the
bordism set is trivial. As a consequence, the mapping class group of $M$ acts trivially on $Emb_{\Z_2}(M,S^6)$. This proves theorem
\ref{invcl}.

{\textsc{Max-Planck-Institut f\"ur Mathematik, Bonn, Vivatsgasse 7, 53111 Bonn, Germany}}\\
E-mail: \nolinkurl{olber@mpim-bonn.mpg.de}

\end{document}